\documentclass[twoside,11pt]{article}

\usepackage{ifthen}
\usepackage{hyperref}
\usepackage{amssymb}
\usepackage{theorem,times,latexsym}
\usepackage{graphics}

\makeatletter
\def\mytag#1{\def\@eqnnum{\normalfont #1 }\global\advance\c@equation\m@ne}
\makeatother

\newtheorem{lemma}{Lemma}

\newtheorem{theorem}{Theorem}
\newtheorem{assumption}{Assumption}
\newtheorem{definition}{Definition}

\newtheorem{example}{Example}
\newtheorem{rem}{Remark}
\newtheorem{proposition}{Proposition}
\newtheorem{cor}[theorem]{Corollary}
\newcommand{\Punkt}{\hspace{-0.5em}{\bf .\ \ }}
\newcommand{\bb}{\mathbb}
\newcommand{\seq}{\longrightarrow}
\newcommand{\imp}{\Rightarrow}
\newcommand{\subs}[3]{{\,}_{#2}#1_{#3}}

\newcommand{\ia}{\alpha}

\newcommand{\id}{\delta}

\newcommand{\qb}{\bar\beta}
\newcommand{\qa}{\bar\alpha}

\newcommand{\qd}{\bar\delta}

\newenvironment{proof}{\vspace{3pt}\noindent
                       \textsc{\em Proof.} }
                       {\hfill$\square$\vspace{3pt}}

\begin{document}

\title{On the existence of Hermitian-harmonic maps
from complete Hermitian to complete Riemannian manifolds
\thanks{Support from the research focus ``Globale Methoden in
der komplexen Geometrie'' under the auspices of ``Deutsche 
Forschungsgemeinschaft'' is gratefully acknowledged}}

\author{
Hans-Christoph 
Grunau\thanks{e-mail:  hans-christoph.grunau@mathematik.uni-magdeburg.de }
\ \ and\ \ 
Marco K\"uhnel\thanks{e-mail: marco.kuehnel@mathematik.uni-magdeburg.de }
}
\maketitle

\begin{abstract}
On non-K\"ahler manifolds the notion of harmonic maps is modified
to that of Hermitian harmonic maps in order to be compatible with the complex 
structure. The resulting semilinear elliptic system is {\it not}
in divergence form. 

The case of noncompact complete preimage and target manifolds is
considered. We give conditions for existence and uniqueness of
Hermitian-harmonic maps and solutions of the corresponding
parabolic system, which observe the non-divergence form of
the underlying equations. Numerous examples illustrate the
theoretical results and the fundamental difference to
harmonic maps.
\end{abstract}

\section{Introduction}
\label{intro}

Let $M$ be a Hermitian manifold of complex dimension
$m$ with Hermitian metric
$
\left( \gamma_{\alpha\, \bar{\beta}} (z) \right)_{\alpha, \beta =1,\ldots,m}
$
and let $N$ be a Riemannian manifold of real dimension $n$
with metric 
$
\left(  g_{j\,k}  (x) \right)_{j,k=1,\ldots,n}
$
and the Levi-Civita-connection, which in local coordinates 
is given by means of  the Christoffel symbols
$\Gamma^j_{k\,\ell} (x)$.
We look for Hermitian harmonic maps $u : M \to N$, which
are defined as solutions of the semilinear elliptic system
\begin{equation}\label{Hermitianharmonicmapsystem}
\gamma^{\alpha \bar{\beta}} \left( \frac{\partial^2 u^{\ell}}{\partial
z^{\alpha} \partial z^{\bar{\beta}}} + \Gamma^{\ell}_{jk}
\frac{\partial u^j}{\partial z^{\alpha}}
\frac{\partial u^k}{\partial z^{\bar{\beta}}}            \right)
=0,\qquad \ell=1,\ldots,n.
\end{equation} 
We focus on the case, where the Hermitian manifold is {\it not}
K\"ahler, and where the system (\ref{Hermitianharmonicmapsystem})
is {\it not} in divergence form. This system was studied first
by Jost and Yau \cite{JostYau}: As they explain, in contrast with
the harmonic map system, this system
is compatible with the holomorphic structure on $M$. They obtain 
beside others existence and uniqueness
results, which cover the Dirichlet problem
for (\ref{Hermitianharmonicmapsystem}) on compact preimage manifolds
with boundary. 
Subsequent work of Chen \cite{Chen} covers the case of {\it target manifolds}
with boundary.
Extensions of existence and uniqueness results for the Dirichlet problem
as obtained in the work of Jost and Yau \cite{JostYau}
to noncompact complete preimage manifolds were
first considered by Lei Ni \cite{LeiNi}. He requires the bilinear
form corresponding to the ``holomorphic Laplace operator'' for {\it functions} 
$u: M \to \mathbb{R}$ 
\begin{equation}\label{holomorphicLaplace}
- \tilde{\Delta } u = 
- 4 \gamma^{\alpha \bar{\beta}} \frac{\partial^2 u}{\partial
z^{\alpha} \partial z^{\bar{\beta}}},
\end{equation}
to be bounded from below by a positive multiple of 
$\int_M u^2 $. 
Such a condition is adequate  -- although very restrictive -- in the
selfadjoint setting, but
 does not really seem to fit in the nonselfadjoint framework
on non-K\"ahler manifolds. 

We impose an invertibility condition on the 
holomorphic Laplace operator between suitably
chosen function spaces, see Assumption~\ref{mainassumption}
below. These function spaces are defined in terms
of decay conditions at ``infinity''. The preimage and
image spaces for the solution operator for 
the holomorphic Laplacian
may be chosen different, and hence our condition
is very flexible and applies to many different situations. 
Even in the selfadjoint setting of harmonic maps this condition
still applies, when $0$ may be a singular value of the Laplace-Beltrami
operator. In this sense, the present note also extends work of Li 
and Tam \cite{LiTam}.

For an extensive discussion
we refer to Subsection~\ref{examples} below. There it will
 become clear that this invertibility condition is even weaker
than assuming that $0$ is not a spectral value for
the holomorphic Laplace operator, defined as a closed unbounded
operator in one fixed function space.

The holomorphic Laplace operator coincides with the usual
Laplace operator if and only if the manifold $M$ is
K\"ahler. That means that we focus on the case, where
the holomorphic Laplacian is {\it not} selfadjoint.

Further we have to assume that there is an initial mapping
 $h:\,M \to N$, such that the Hermitian-harmonic differential
operator, applied to $h$, decays suitably at $\infty$. Then we can 
show the existence of a Hermitian-harmonic map $u:M \to N$, which
is homotopic to $h$ and which approaches $h$ at $\infty$. This main
result is contained in Subsection~\ref{existence}

In \cite{LiTam} examples of harmonic diffeomorphisms are given 
which are homotopic to the identity. One might expect
to see similar examples for Hermitian harmonic maps. 
In Section \ref{Limits} we prove that in a series of manifolds including the ones
in \cite{LiTam}, it is not possible for the identity to satisfy the decay condition mentioned above.
We believe that in those cases there do not exist Hermitian 
harmonic diffeomorphisms homotopic to the identity. It is not only in
this respect that the complex structure of the preimage manifold and
the nonselfadjoint principal part of the elliptic system 
complicate the construction of relevant examples.

Originally, existence of harmonic as well as of Hermitian-harmonic 
maps was proved via the seeming detour of the corresponding parabolic
equations. The reason is the lack of compactness properties of the
underlying elliptic systems. That this approach works out also
for non divergence form systems with a nonlinearity 
quadratic in the gradient, was observed first in \cite{vonwahl}.
In \cite{JostYau}, the parabolic method was applied to the study
of Hermitian-harmonic maps, and the required stability and
convergence properties in $C^0$-norms were found.
In the present paper, as well as in \cite{LeiNi}, the exhaustion procedure
will work directly on the elliptic level. Nevertheless it is interesting
to know, whether solutions to (\ref{Hermitianharmonicmapsystem}) 
may be obtained as limits for $t \to \infty$ of the corresponding
parabolic system also in our noncompact situation. This
question is addressed and answered in Section~\ref{heatequation}.
To ensure convergence we need to impose a decay condition on the
linear heat operator
$$
\left( \frac{\partial}{\partial t} -\tilde{\Delta}\right)u,
$$
which is related to the invertibility condition for the 
holomorphic Laplace operator.
This decay condition is discussed and illustrated 
in Subsection~\ref{parabolicexamples} with help 
of the same series of examples as for the elliptic system.

\section{The elliptic Hermitian-harmonic map system}
\label{ellipticsystem}

\subsection{Preliminaries}
\label{preliminary}

In this section, after explaining the notation,
we collect some basic results from the
fundamental papers on Hermitian-harmonic maps by 
Jost and Yau \cite{JostYau} and Lei Ni \cite{LeiNi}.

First we specify and explain our notation.
Let $M$ be a complete Hermitian manifold of complex dimension
$m$ with Hermitian metric
$$
\left( \gamma_{\alpha\, \bar{\beta}} (z) \right)_{\alpha, \beta =1,\ldots,m}
$$
in local coordinates.
By $ \gamma^{\alpha\, \bar{\beta}}$ we denote the transposed inverse
matrix
$$
\sum_{\sigma=1,\ldots,m}
\left(\gamma^{\alpha\, \bar{\sigma}}
     \gamma_{\beta\, \bar{\sigma}} (z) \right)
       = \delta^{\alpha}_{\beta}.
$$
With respect to this metric,  the length of
a ho\-lo\-mor\-phic tan\-gen\-tial vector $w = \left(w^1,\ldots,w^m \right)$ 
at $z \in M$ in local coordinates is given by
$$
\| w \|^2 = \sum_{\alpha , \beta =1,\ldots,m}
      w^{\alpha }  \gamma_{\alpha\, \bar{\beta}} (z) \bar{w}^{\beta}  .
$$
Furthermore, let $N$ be a complete Riemannian manifold of real dimension $n$
with metric 
$$
\left(  g_{j\,k}   \right)_{j,k=1,\ldots,n}
$$
in local coordinates, its inverse
$$
\sum_{\ell = 1, \ldots, n} g_{j\, \ell } g^{\ell\, k} = 
  \delta_j^{k}
$$
and the Christoffel symbols
$$
\Gamma^j_{k\,\ell} = \frac{1}{2} \sum_{s=1}^n g^{j s}
   \left( \frac{\partial g_{\ell s  } }{\partial x^k }
         + \frac{\partial g_{ s k } }{\partial x^{\ell} }
         - \frac{\partial g_{ k \ell } }{\partial x^s }
        \right).
$$
While on the target manifold $N$, we consider the Levi-Civita
connection of the metric, we choose a different connection
on the preimage manifold $M$. We choose a 
suitable holomorphic torsion
free connection such that the ``holomorphic Laplace operator''
takes the form as above in (\ref{holomorphicLaplace}).

Further we need to define the tension field for any smooth map
$u:M \to N$ according to the chosen connections
\begin{equation}\label{tensionfield}
\left( \sigma (u) \right)^\ell \, := \,
\gamma^{\alpha \bar{\beta}} \left( \frac{\partial^2 u^{\ell}}{\partial
z^{\alpha} \partial z^{\bar{\beta}}} + \Gamma^{\ell}_{jk}
\frac{\partial u^j}{\partial z^{\alpha}}
\frac{\partial u^k}{\partial z^{\bar{\beta}}}            \right),
\qquad \ell=1,\ldots,n.
\end{equation}

The first result we need to mention
concerns the energy density function $e(u)$, 
which for any smooth map $u:M \to N$ is defined in local
coordinates as follows
\begin{equation}\label{energydensity}
e(u) := \left( g_{jk} \circ u \right) \gamma^{\alpha \bar{\beta}} 
\frac{\partial u^j}{\partial z^{\alpha}}
\frac{\partial u^k}{\partial z^{\bar{\beta}}}
\end{equation}
If we assume $u:M \to N$ to be a Hermitian harmonic map and $N$ 
to have nonpositive sectional curvature, 
then according to \cite[p. 225, formula (5)]{JostYau}, 
for any relatively compact open set 
$\Omega \subset M$ we have 
the following differential inequality
\begin{equation}\label{Lemma3.1}
-\tilde{\Delta} e (u) \le  
     C (\Omega)  e(u).
\end{equation}
The constant $C$ is expected to blow up in general, when
$\Omega$ is approaching  $M$. For the reader's convenience we sketch the proof 
of (\ref{Lemma3.1}) in Appendix~\ref{appendixb}.

One should observe that by the Hopf-Rinow-theorem 
(see e.g. \cite[1.37]{Aubin}) the compact subsets of $M$ 
are precisely the bounded closed sets.

The next important result is due to Lei Ni 
\cite[Corollary 3.5]{LeiNi}. For this we need first to
explain the geodesic homotopy distance between 
two smooth homotopic maps $u$ and $v:M \to N$. 
Let us recall a result of von Mangoldt-Hadamard-Cartan.
Fix a homotopy $H$
between $u$ and $v$, then, since the target manifold is nonpositively
curved,
for any $z \in M$ there is precisely one geodesic arc connecting
$u(z)$ and $v(z)$ in the same homotopy class as the original arc 
given by $H$.  Moreover this geodesic arc is length minimizing.
See e.g. \cite[Lemma 8.7.1]{Jost2}.
The  geodesic homotopy distance 
$$
\rho := \rho (z) := \rho (u(z), v(z))
$$
is defined as the length of this geodesic arc.

According to \cite[Corollary 3.5]{LeiNi},
$\rho$ satisfies the following fundamental 
differential inequality holds true:
\begin{equation}\label{Corollary3.5}
- \tilde{\Delta} \rho \le 
     4 \left( \| \sigma (u)\| + \| \sigma (v) \| \right).
\end{equation}
In the next section, we will construct Hermitian-harmonic maps
by an exhaustion procedure and by solving a boundary value problem for
(\ref{Hermitianharmonicmapsystem}) on compact submanifolds of $M$. 
The above estimate will
turn out to be essential for getting first estimates for the approximate
solutions to (\ref{Hermitianharmonicmapsystem}).

%
%
%
%

\subsection{Existence and uniqueness results}
\label{existence}

We first introduce spaces of suitably decaying functions (at
``infinity''), which are adequate in our nonselfadjoint 
and noncompact framework.

\begin{definition}\Punkt \ 
For $\mu > 0$, let
\begin{eqnarray}
C^0_{\mu} (M) & := & \left\{f:M \to \mathbb{R};\ f \mbox{\ is continuous
   and }
   \right.\\
&& \left. \mbox{\ 
 there
exists $z_0 \in M$ and a constant $C=C(f)$ such that\ }
|f(z)| \le C\, \left( 1+ d(z,z_0)\right)^{-\mu} \right\}. \nonumber
\end{eqnarray}
\end{definition}

\begin{assumption}[Invertibility  of the holomorphic
Laplace operator]\Punkt\label{mainassumption}
\ \ \ \\
We assume that there exist positive numbers $\mu,\mu'>0$ such that
for every $f\in C^0_{\mu} (M)$, there exists precisely one solution
$u \in C^0_{\mu'} (M)$ of
$$
- \tilde{\Delta} u = f \mbox{\ in\ } M.
$$
\end{assumption}

\begin{theorem}[Existence and uniqueness of Hermitian harmonic maps]\Punkt\label{existenceresult}
Assume that $M$ is a noncompact complete Hermitian manifold such that for the
holomorphic Laplace operator $- \tilde{\Delta} $ on $M$, the 
Assumption~\ref{mainassumption} is satisfied with positive numbers
 $\mu,\mu'>0$. Further let $N$ be a 
complete Riemannian manifold with nonpositive sectional curvature
and $h:M \to N$ a smooth map with $\| \sigma (h) \|  \in  C^0_{\mu} (M)$.

Then there exists a Hermitian harmonic map $u:M \to N$, which is homotopic
to $h$. Moreover, if $\rho$ denotes the homotopy distance between $u$ and $h$,
we have $\rho \in  C^0_{\mu'} (M)$. Finally, in this class, the solution
is unique.
\end{theorem}

\vspace{4mm}
\noindent
{\it Proof.} The fundamental idea is as in the paper \cite{LeiNi}.
Here, however, we replace the ``selfadjoint'' tools by the appropriate
nonselfadjoint analogues. Let $\left( \Omega_k \right)_{k \in \mathbb{N}}$ 
be a compact smooth exhaustion of $M$.
According to Theorem 6 of the paper \cite{JostYau} by J. Jost and
S.-T. Yau, there exist solutions $u_k: \Omega_k \to N$ of the 
Dirichlet problems
\begin{equation}\label{compactDirichletproblems}
\left\{ \begin{array}{ll}
\displaystyle \sigma(u_k) = 0  &\mbox{\ in\ } \Omega_k,\\[2mm]
\displaystyle u_k = h           &\mbox{\ on\ } \partial\Omega_k,\\[2mm]
\displaystyle u_k \mbox{\ homotopic to\ } h, &\mbox{\ with respect to\ }
 \partial\Omega_k.
\end{array} \right.
\end{equation}
In order to show convergence of $(u_k)$ to a Hermitian harmonic map
$u:M \to N$, it is enough to prove local boundedness of the 
energy density functions $e(u_k)$. As in \cite{LeiNi} we start
with global bounds for the homotopy distances $\rho_k$ between
$u_k$ and $h$ and  $\rho_{k,\ell}$ between $u_k$ and $u_{\ell}$.

We first introduce a comparison function, the existence of which
is ensured by Assumption~\ref{mainassumption}: Since $\| \sigma (h) \|$
is assumed to be in $ C^0_{\mu} (M)$, we find a smooth function
$V \in  C^0_{\mu'} (M)$, such that
\begin{equation}\label{comparisonfunction}
- \tilde{\Delta}  V = 4 \| \sigma (h) \| \mbox{\ in\ } M.
\end{equation}
In particular, $V(z)$ decays uniformly, as $d(z,z_0) \to \infty$. Together
with the strong maximum principle, which can be easily proven by 
passing to local coordinates and by exploiting the connectedness of $M$,
this gives first $V \ge 0$ and then by repeating the argument:
\begin{equation}\label{positivity}
V > 0.
\end{equation}

By (\ref{Corollary3.5}), the coincidence of $\rho_k$ and $h$ on
$\partial \Omega_k$  and (\ref{positivity}), we find the following
inequalities for the homotopy distance $\rho_k$ between
$u_k $ and $h$:
\begin{eqnarray*}
- \tilde{\Delta} \rho_k & \le & 4 \| \sigma (h) \| =
   - \tilde{\Delta}  V  \mbox {\ in\ } \Omega_k,\\
\rho_k|\partial \Omega_k &=& 0 < V|\partial \Omega_k.
\end{eqnarray*}
{From} the maximum principle, we get the uniform bound:
\begin{equation}\label{diffineqhomotopydist}
0 \le \rho_k \le V,
\end{equation}
where $V \in  C^0_{\mu'} (M)$ is the comparison function,
introduced in (\ref{comparisonfunction}) above.

In a second step we will exploit the differential inequality 
(\ref{Lemma3.1}) for the energy density $e(u_k)$ of the 
approximate Hermitian harmonic maps $u_k$.

We take a local $L^1$-bound for $e(u_k)$ from
\cite[pp. 344/345]{LeiNi}:
For some fixed $z_0 \in M$ and any $R >0$, we have with a suitable constant
\begin{equation}\label{L1bound}
\int_{B_R (z_0)} e (u_k) \le C.
\end{equation}
This bound holds true also  in our situation as we have shown
the maximum bound (\ref{diffineqhomotopydist}) for $\rho_k$
above.

Eventually from this local $L^1$-bound (\ref{L1bound}), 
we get local $L^{\infty}$-bounds
by making use of the  local maximum principle 
\cite[Theorem 9.20]{GilbargTrudinger} for elliptic operators,
which are not in divergence form.
First we work in sufficiently small open sets of $M$, where simply one chart is sufficient. 
The holomorphic Laplace operator in these local coordinates
satisfies the assumptions of the local maximum principle and we exploit
the differential inequality
$$
- \tilde{\Delta} e (u_k) \le C_{\mbox{\footnotesize loc}} e(u_k).
$$
See (\ref{Lemma3.1}); the constant can be found at least
on compact subsets of $M$. Second, as
by the Hopf-Rinow theorem (see e.g. \cite[1.37]{Aubin}), all the 
$\overline{\Omega}_{\ell}$ are compact, we get there with help
of a bootstrapping argument uniform $C^{2,\alpha}$-bounds and hence 
convergence to a smooth solution $u$ of the
Hermitian harmonic map system (\ref{Hermitianharmonicmapsystem}).

It is easy to see that $u$ and $h$ are homotopic. For this purpose
we extend $u_k:\Omega_k \to N$ by $h$ to a continuous mapping
$\tilde{u}_k: M \to N$. Further let $\tilde{u}_0:=h$. 
Obviously, $\tilde{u}_k$ and $\tilde{u}_{k+1}$ 
are homotopic; for $k \in \mathbb{N}_0$ let 
$H_k: \left[\frac{1}{k+2},\frac{1}{k+1}\right] \times M \to N$
be continuous with $H_k \left( \frac{1}{k+1} , \, . \,\right) = \tilde{u}_k$
and $H_k \left( \frac{1}{k+2}, \, . \, \right) = \tilde{u}_{k+1}$.
Defining
\begin{eqnarray*}
H &:& [0,1] \times M \to N,\\
H \left( t, \, . \, \right) &=& \left\{ \begin{array}{ll}
\displaystyle
H_k \left( t, \, . \, \right), &\mbox{\ if\ } t 
         \in \left[\frac{1}{k+2},\frac{1}{k+1}\right],\\[2mm]
\displaystyle         
u, &\mbox{\ if\ } t=0,
\end{array}\right. 
\end{eqnarray*}
we get a homotopy between $u$ and $h$.

We conclude from (\ref{diffineqhomotopydist}) and locally 
uniform convergence that $0 \le \rho \le V$ and hence
$\rho \in C^0_{\mu'} (M)$.\hspace*{4mm}

Finally we prove uniqueness of the solution $u$ with the 
mentioned properties. Let $\tilde u:M\seq N$ be an arbitrary Hermitian-harmonic
map of class $C^0_\mu(M)$ homotopic to $h$, such that 
$\rho(\tilde u,h)\in C^0_{\mu'} (M)$. By (\ref{Corollary3.5}) we know
$$-\tilde\Delta\rho(u,\tilde u)\le 0$$
and, by the previous arguments, that
$$0\le\rho(u,\tilde u)\le\rho(u,h)+\rho(\tilde u,h)\in C^0_{\mu'}(M).$$
Hence for every $\varepsilon>0$ outside a sufficiently large 
ball $B_R(z_0)$ around an arbitrary $z_0\in M$ we have
$$\rho(u,\tilde u)\le \varepsilon.$$
By the maximum principle this implies $\rho(u,\tilde u)\le
\varepsilon$ on all of $M$ for every $\varepsilon>0$
and hence $\rho(u,\tilde u)=0$.
This implies $u=\tilde u$.
 \ \ \ \ \ \ \ \ \ \ \ \ \ \ \ \ \   \hfill  $\square$

\subsection{Examples}
\label{examples}

First, with help of some examples, 
 we want to discuss the invertibility condition on 
 the holomorphic Laplace operator, i.e. Assumption~\ref{mainassumption}.
We are basing our first examples on the following simple result:

\begin{lemma}\Punkt\label{basic}
Let $n >2$, $\alpha \in \left( 0,\frac{n}{2} -1 \right)$. Then,
for every continuous $f :\mathbb{R}^n \to \mathbb{R}$ with 
$|f(x)| \le C\, (1+|x|^2)^{-\alpha-1}$, we find precisely one
strong solution $u:\mathbb{R}^n \to \mathbb{R}$ of
\begin{equation}\label{entiresolution}
- \Delta u = f \mbox{\ in\ } \mathbb{R}^n,
\end{equation}
such that
$$
|u(x)| \le C\, (1+|x|^2)^{-\alpha}.
$$
\end{lemma}

\vspace{4mm}
\noindent
{\it Proof.} We define
$$
v(x) := (1+|x|^2)^{- \alpha}
$$
as a barrier function and calculate:
\begin{eqnarray*}
- \Delta v (x) &=& 2 \alpha n (1 + |x|^2)^{-\alpha -1}
                  - 4 \alpha (\alpha + 1) |x|^2 (1 + |x|^2)^{-\alpha -2}\\
     & \ge & c_{n,\alpha} (1 + |x|^2)^{-\alpha -1};
\end{eqnarray*}
where the positive constant $c_{n,\alpha}$ is given by
$$
c_{n,\alpha} = 4 \alpha \left( \frac{n}{2} -(\alpha + 1)\right).
$$
In order to find a solution to (\ref{entiresolution}), we solve
the corresponding Dirichlet problems with homogeneous boundary data
on the balls $B_k$ around the origin with radius $k$. As a suitable
multiple of $v$ will serve as a barrier function for the 
approximate solutions $|u_k|$, after selecting a suitable
subsequence we will have local convergence in
$ C^0$ and weakly in $W^{2,p}$ for arbitrarily large $p$ against an entire solution
of (\ref{entiresolution}), obeying the same bound $C\, v(x)$.

Uniqueness is immediate from Liouville's theorem. 
\hfill $\square$

\begin{example}\Punkt \ 
Let us consider $M=\mathbb{C}^m$, $m \ge 2$, with the
standard euclidean metric, so that the holomorphic Laplacian is
also the standard one: $\Delta_e$. Then, according to the previous lemma,
Assumption~\ref{mainassumption} is satisfied with any 
$\mu \in (2,2m)$ and $\mu'=\mu-2$.
\end{example}

\vspace{4mm}\noindent
In this example, the holomorphic Laplace operator
is selfadjoint. Although we do not focus on this case here,
this observation shows: Even if $0$ is a singular value of the 
Laplace operator, our invertibility assumption may still be satisfied.

In order to cover also nonselfadjoint examples,
we would like to equip $M=\mathbb{C}^m$, $m \ge 2$ with the
conformal metric
$$
\gamma_{\alpha\, \bar{\beta}}(z) = 
            (1+|z|^2)^{-1}\delta_{\alpha\,\beta}.
$$
The holomorphic Laplace operator then becomes
$$
-  \tilde{\Delta} = - (1+|z|^2) \Delta_e
$$
with $\Delta_e$ being the euclidean Laplace operator. In 
$L^2 \left(\mathbb{C}^m,\gamma \left( \frac{i}{2}\right)^m 
  (d\, z_1 \wedge d \, \bar{z}_1) \wedge \ldots \wedge  (d\, z_m \wedge d \, \bar{z}_m) \right)
= L^2 \left( \mathbb{C}^m,\left( (1+|z|^2)^{-m} \left( \frac{i}{2}\right)^m 
  (d\, z_1 \wedge d \, \bar{z}_1) \wedge \ldots \wedge  (d\, z_m \wedge d \, \bar{z}_m) \right)\right)$,
the holomorphic Laplacian is not selfadjoint.

Since
$$\omega=\frac{i}{2}(1+|z|^2)^{-1}\sum dz_{\alpha}\wedge d\bar{z}_{\alpha}$$ we compute
$$d\omega=\frac{i}{2} (1+|z|^2)^{-2}
   \sum_{\alpha, \beta}2(\bar{z_\beta}dz_\beta\wedge dz_\alpha\wedge d\bar{z_\alpha}+z_\beta dz_\alpha\wedge d\bar{z_\alpha}\wedge
d\bar{z_\beta})\not=0.$$
This means, that $(M,\gamma)$ is not a K\"ahler manifold, what is important, since otherwise 
Hermitian-harmonic maps are harmonic.

Again, Lemma~\ref{basic} shows, that for any smooth $f$ with
$|f(z)| \le C (1+|z|^2)^{-\alpha}$, $\alpha \in (0,m-1)$, we find a solution $u$
of
$$
- \tilde{\Delta} u = f \mbox{\ \ \  in\ } \mathbb{C}^m
$$
with $|u(z)| \le C (1+|z|^2)^{-\alpha}$. However, with this metric
we have
$$
d(z,0) \sim \log(1+|z|^2),\qquad |z| \sim \exp (d(z,0)) -1.
$$
This example doesn't fall under our formulation of
Assumption~\ref{mainassumption}. However it shows that the choice
of the metric
$$
\gamma_{\alpha\, \bar{\beta}}(z) = 
            (1+|z|^2)^{-1}\delta_{\alpha\,\beta}.
$$
and of the corresponding holomorphic Laplace operator 
$$
-  \tilde{\Delta} = - (1+|z|^2) \Delta_e
$$
may be reasonable. Since $\log (1+|z|^2) \sim d(0,z)$, where $|z|$ is the 
euclidean norm and $d(z,0)$ the distance in our metric to the origin, we should find 
a refinement of Lemma~\ref{basic}, which involves logarithmic terms:

\begin{lemma}\Punkt\label{notsobasic}
Let the dimension be $n >2$ and let $\alpha > 0$ be a real
number. Then 
for every $f :\mathbb{R}^n \to \mathbb{R}$ with 
$|f(x)| \le C\, \left( \log(2+|x|^2) \right)^{-\alpha-1}$, we find precisely one
solution $u:\mathbb{R}^n \to \mathbb{R}$ of
\begin{equation}\label{entirelogarithmicsolution}
- (1+|x|^2) \Delta u = f \mbox{\ in\ } \mathbb{R}^n,
\end{equation}
such that
$$
|u(x)| \le C\, \left( \log (2+|x|^2) \right)^{-\alpha}.
$$
\end{lemma}

\vspace{4mm}
\noindent
{\it Proof.} Similarly as in the proof of Lemma~\ref{basic} we look
for a suitable comparison function. First let us work with an
auxiliary number $A \ge 2$, which will be fixed in the course 
of the following calculations. We define 
$$
v(x):= \left(  \log (A + |x|^2 )\right)^{- \alpha}
$$
and calculate:
\begin{eqnarray*}
- \Delta v (x) &=& 2 \alpha \left( \log (A + |x|^2)\right)^{- \alpha -1} 
     \left( \frac{n}{A+|x|^2} - 2\frac{|x|^2}{(A + |x|^2 )^2}\right)\\
    &&  - 4 \alpha (\alpha +1) \left( \log (A + |x|^2)\right)^{- \alpha -2}
       \, \frac{|x|^2}{(A+ |x|^2)^2} \\
   & \ge& 2 \alpha (n -2 )  \left( \log (A + |x|^2)\right)^{- \alpha -1} 
                 \frac{1}{A + |x|^2}
         \left\{  1- 2\frac{\alpha+1}{n-2}  
      \left( \log (A + |x|^2)\right)^{ -1}  \right\}\\
    & \ge & \frac{\alpha (n -2 )  }{A}
      \left( \log (A + |x|^2)\right)^{- \alpha -1} 
                 \frac{1}{1+ |x|^2},
\end{eqnarray*}
provided $A$ is chosen large enough in dependence on $\alpha >0$
and $n > 2$. As in the proof of Lemma~\ref{basic}, we have now:
For every continuous function $f$ with $|f(x)| \le C 
\left( \log (A + |x|^2)\right)^{- \alpha -1} $ we have precisely
one solution $u$ of
$
- (1+ |x|^2) \Delta u (x) = f(x) \mbox{\ in\ }\mathbb{R}^n
$
with $| u(x) | \le C \left( \log (A + |x|^2)\right)^{- \alpha } $.
But since the strictly positive function 
$ (0,\infty) \ni r \mapsto \log (2+ r^2)/\log(A+ r^2)$
is bounded from above and below, this immediately gives 
the statement of our lemma. \hfill $\square$

\begin{example}\Punkt \ \ \label{basicex}
Let $M=\mathbb{C}^m$, $m \ge 2$ be equipped with the
conformal metric
$$
\gamma_{\alpha\, \bar{\beta}}(z) = 
            (1+|z|^2)^{-1}\delta_{\alpha\,\beta},
$$
such that the holomorphic Laplace operator is
$$
-  \tilde{\Delta} = - (1+|z|^2) \Delta_e
$$
with $\Delta_e$ being the euclidean Laplace operator. 
Then $ -  \tilde{\Delta} $ satisfies the invertibility
condition Assumption~\ref{mainassumption} with any
$\mu > 1$ and $\mu' = \mu -1$. 
\end{example}

                                                                                                                                                                              \vspace{4mm}\noindent
The second purpose of this subsection is to discuss the decay condition
on $\| \sigma (h) \| \in C^0_{\mu} (M)$. For this we construct some
prototype manifolds $M$ and $N$ and  suitable ``initial maps''
$h: M \to N$.

\begin{example}\Punkt \
{\rm 
On ${\bb R}^2$ 
the rotational symmetric metric $g_0=dr^2+(r^2+r^4)d\phi^2$ has strictly
negative curvature. If we now choose $N={\bb R}^2\times{\bb R}^2$ with the metric $g=pr_1^*g_0+pr_2^*g_0$. where
$pr_i:N\seq{\bb R}^2$ denotes the projections onto the 
$i$-th copy of ${\bb R}^2$, then $(N,g)$ has nonpositive sectional
curvature. 

As for the manifold $M$ we first choose $\tilde M={\bb C}^2$ with the Hermitian metric
$$\tilde\gamma=\frac{1}{1+|z|^2}(dz_1\otimes d\bar z_1+dz_2\otimes d\bar z_2).$$
Then it is easy to see that the geodesic length $d(z,0)\sim \log(1+|z|^2).$

Now $M:=\tilde M\setminus B_1(0)$ shall be regarded as a manifold with boundary $\partial B_1(0)$. The proof
of the theorem works in the same way for this $M$ where additionally $u=h$ on $\partial B_1(0)$ can be
satisfied. 

If we define $h:M\seq N$ via 
$$h(z)=\frac{z}{1+|z|^2},$$
the norm of the tension field $\|\sigma(h)\|$ can be computed to be
$$\|\sigma(h)\|=\frac{|z|(7+2|z|^2)}{2(1+|z|^2)^2}\le\frac{7|z|}{2(1+|z|^2)},$$
and hence $\|\sigma(h)\|\in C^0_\mu(M)$ for every $\mu>0$.

Applying Example \ref{basicex} for $\mu>1$ yields by Theorem \ref{existenceresult} a Hermitian-harmonic
map $u:M\seq N$ homotopic to $h$ with $u=h$ at $\partial B_1(0)$ and 
approaching $0$ at infinity. In particular, $u$ is
not a constant map.}
\end{example}

%

\begin{example}\label{zweienden}\Punkt \
{\rm
Let $M=({\bb S}^1)^{2m-1}\times(-1,1)$ equipped with the following complex structure: Denote $H:=\{z\in{\bb C}|\,|\Im(z)|<1\}$ and take
$({\bb C}^{m-1}\times H)/\Gamma\cong M$, where $\Gamma$ is the cartesian lattice of rank $2m-1$. 

For the choice of the metric, denote by $s$ the noncompact parameter with range
$s\in(-1,1)$. Then the metric
$$\tilde\gamma_{\alpha\bar\beta}:=f(s)\delta_{\alpha\bar\beta}$$
in cartesian coordinates 
is $\Gamma$-invariant and not K\"ahler unless $f$ is constant. 
We denote by $\gamma$ the induced metric on $M$. We choose $\delta>0$
and
$$
f(s):=\left\{\begin{array}{ll}
\delta^2(1-|s|)^{-2\delta-2}&\mbox{\ for\ }1/2<|s|<1,\\
a(s)&\mbox{\ for\ }|s|\le 1/2,
\end{array}\right. 
$$
such that $a(s)>0$ for all $|s|\le 1/2$ and $f\in C^\infty(-1,1)$. 

Since for $|s|$ close to $1$ one has
$d(z,0)\sim (1-|s|)^{-\delta}-1:=\tilde d(s)$, the metric $\gamma$ is complete.

We can prove that $-\tilde\Delta(1+\tilde d)^{-{\mu'}}>C(1+\tilde d)^{-{\mu'}-2}$ as long as $\delta{\mu'}<1$ and $|s|>1/2$: Since
$$(1+\tilde d(s))^{-{\mu'}}=(1-|s|)^{\delta{\mu'}},$$
we get even 
\begin{eqnarray*}-\tilde\Delta(1+\tilde d(s))^{-{\mu'}}&=&-\frac 1{\delta^2}(1-|s|)^{2\delta+2}\frac{\partial^2}{\partial s^2}(1-|s|)^{\delta{\mu'}}\\
&=&\frac{{\mu'}(1-\delta{\mu'})}{\delta}(1-|s|)^{2\delta+2}(1-|s|)^{\delta{\mu'}-2}\\
&=&\frac{{\mu'}(1-\delta{\mu'})}{\delta}(1-|s|)^{\delta({\mu'}+2)}\\
&=&\frac{{\mu'}(1-\delta{\mu'})}{\delta}(1+\tilde d(s))^{-{\mu'}-2}.\end{eqnarray*}
 
Now we remark that $b(s)=1+\varepsilon (1/4-s^2)$ satisfies $-\tilde\Delta b>0$ 
for every $\varepsilon>0$. We define
$$
v(s)=\left\{\begin{array}{ll}
(1+\tilde d(s))^{-{\mu'}}&\mbox{ for }|s|>1/2,\\
(1+\tilde d(1/2))^{-{\mu'}}b(s)&\mbox{ for }|s|\le 1/2.
\end{array}\right.$$
Then we compute for $\phi\in C^\infty_0(M), \phi\ge 0$
\begin{eqnarray*}
\lefteqn{
\int\left(-\tilde{\Delta}^*\phi\right)v  \, f^m\, dx \ge
\left(\phi f^{m-1} \right)(1/2)\left(\frac{\partial
v}{\partial s}\left(\frac 12-0\right)- \frac{\partial v}{\partial s}\left(\frac
12+0\right)\right)}\\
&&+\left(\phi f^{m-1} \right)(-1/2)\left(\frac{\partial
v}{\partial s}\left(-\frac 12-0\right)- \frac{\partial v}{\partial
s}\left(-\frac 12+0\right)\right)\\
&=&\left( \left( \phi f^{m-1}\right)(1/2)+  \left( \phi f^{m-1}\right)
(-1/2) \right) \left(-\varepsilon+{\mu'}\frac{\tilde
d'\left(1/2\right)}{1+\tilde d\left(1/2\right)}\right) \left(1+\tilde
d\left(1/2\right)\right)^{-{\mu'}}\\ &\ge&0
\end{eqnarray*}
for $\varepsilon$ sufficiently small. 
Hence $v(s)\in C^0_{\mu'}(M)$ is a supersolution and we proceed like 
before to prove the validity of Assumption \ref{mainassumption}
with $\mu'=\mu-2$, $\mu>2$ and $\delta (\mu - 2)<1$. Hence we have proved:

\begin{lemma}\Punkt\label{Zweienden}Let $M:=(\bb S^1)^{2m-1}\times(-1,1)$ be 
like in Example \ref{zweienden}. If $\delta>0,\mu>2$
and $\delta(\mu-2)<1$, then Assumption \ref{mainassumption} is valid with $\mu':=\mu-2$.
\end{lemma} 

\vspace{4mm}
\noindent
Now we construct a starting map $h$. The idea is to fix the values in both infinite edges and to 
interpolate such that
$\|\sigma(h)\|\in C^0_\mu(M)$.

For this purpose denote $N=B_1(0)\subset{\bb R}^n$ equipped with the Poincar\'e metric 
$g=\frac{1}{(1-r^2)^2}\delta_{\alpha\bar\beta}$.}

\end{example} 

\begin{proposition}\Punkt If $M$ is like in Example \ref{zweienden}, 
$N$ is the unit ball with the Poincar\'e metric, 
and $\tilde h:{\bb C}^{m-1}\times H\seq{\bb R}^n$ a $\Gamma$-invariant 
$C^2$-map with bounded first and second
derivatives and the image $\tilde h({\bb C}^{m-1}\times H)\subset N$
being precompact in $N$, then there is a Hermitian-harmonic map $u:M\seq N$ 
homotopic to the quotient map $h:M\seq N$.
\end{proposition}

\begin{proof}It suffices to prove that $\|\sigma(h)\|\in C^0_\mu(M)$ for some $\mu>2$. 
For this purpose we choose 
$2<\mu<2+\frac 1\delta$. Then Lemma \ref{Zweienden} shows that Assumption \ref{mainassumption} is valid.
We note that $|\Gamma^l_{jk}|\le\frac{r}{1-r^2}\le\frac{1}{1-r^2}$ for the given Poincar\'e metric, 
which is an easy calculation.
By assumption, 
$$\left|\frac{\partial^2}{\partial z^\alpha\partial z^{\bar\beta}}\tilde h^j\right|\le C_1,
        \mbox{  }\left|\frac{\partial}{\partial z^\alpha}\tilde h^j\right|\le C_1,$$
and $r(x):=|\tilde h(x)|\le q<1$ in ${\bb C}^{m-1}\times H$.
Now we can estimate
\begin{eqnarray*}
\|\sigma(h)\|^2&\le& C_2\frac{(1-|s|)^{4+4\delta}}{(1-r^2)^2}\left(C_1+\frac{C_3}{(1-r^2)}\right)^2\\
&\le &C_4(1-|s|)^{4+4\delta}\\
&\le &C_4(1-|s|)^{2\delta\mu}\\
&=&C_4(1+\tilde d(s))^{-2\mu},
\end{eqnarray*}
if $2<\mu\le 2+\frac 2\delta$. By our choice even $2<\mu<2+\frac 1\delta$ holds.
\end{proof}

\subsection{Negative Results}\label{Limits}

Harmonic maps sometimes may be thought of as diffeomorphisms or deformations of the identity in an appropriate 
setting. If, for example, $M=N$ is the unit ball equipped with the Poincar\'e metric, in \cite{LiTam} it is proved that there is a
harmonic map $u:M\seq N$ homotopic to the identity. This suggests to choose $h$ as an identity map and to use Theorem 
\ref{existenceresult} in order to obtain a Hermitian-harmonic map homotopic to $h$. 
This idea fails in many examples, in particular,
if $M=N$ is the unit ball with the Poincar\'e metric. 
We will prove that in this case the assumptions of Theorem
\ref{existenceresult} are not satisfied. We have to leave open  whether there
are Hermitian-harmonic maps homotopic to the identity.

Since we are now going to inquire into rotational symmetric metrics, let us collect some basic knowledge. 

\begin{lemma}\Punkt\label{lines}Let $B_r(0)\subset{\bb R}^k$
be equipped with a rotational symmetric Riemannian metric $\gamma$. 
Let $x\in B_r(0)$ and $\Gamma$ be a geodesic connecting 
$x$ and $0$. Then $\Gamma$ is a line segment.
\end{lemma}

\begin{proof}First we note that by rotational symmetry the geodesic equations tell us that the
line segment between $x$ and $0$ is a geodesic. Now take $y\in\Gamma$ near $0$ such that there is only one geodesic
through $y$ and $0$. This has to be the line segment connecting $y$ and $0$. Since the line through $y$ and $0$
is the unique geodesic with tangent direction $\Gamma'(y)$ in $y$, we conclude that $\Gamma$ is the line segment between
$x$ and $0$.  
\end{proof}

\vspace{4mm}
\noindent
If $0\notin\Omega$ we obtain a somewhat weaker result:

\begin{lemma}\Punkt\label{block}Let $I\subset{\bb R}^+$ be an open interval and $I\times S^{k-1}\cong\Omega\subset{\bb R}^k$
be an annulus equipped with a rotational symmetric Riemannian metric $\gamma$ of the form
$\gamma={\mbox pr}_1^*\gamma_r+{\mbox pr}_2^*\gamma_\phi$ ('polar block form'). 
Let $x,y\in M$ be collinear with $0$. Then the shortest geodesic between $x$ and $y$ is a line segment.
\end{lemma}
 
\begin{proof}
By assumption,
$$\gamma=a(r)dr^2+b_{ij}(r)d\phi_id\phi_j,$$
with $a>0,(b_{ij})>0$. If $\Gamma:[0,1]\seq\Omega$ is a path connecting $x$ and $y$, then
$$l(\Gamma)=\int_0^1\sqrt{a\left(\frac{d\Gamma_r}{ds}\right)^2+b_{ij}\frac{d\Gamma_{\phi_i}}{ds}\frac{d\Gamma_{\phi_j}}{ds}}ds\ge
\int_0^1\sqrt{a\left(\frac{d\Gamma_r}{ds}\right)^2}ds=l(L),$$
if $L$ denotes the line segment between $x$ and $y$.
\end{proof}

\begin{rem}\Punkt Note that the polar block form condition of Lemma \ref{block} is satisfied if $\gamma$ is conformal to the euclidean metric.
\end{rem}

\vspace{4mm}
\noindent
First we show the positive result that the Poincar\'e ball is within the range of Assumption \ref{mainassumption}:

\begin{example}\Punkt Let $M=D^4:=\{z\in{\bb C}^2|\quad|z|<1\}$ equipped with the Poincar\'e metric
$\gamma:=\frac{4}{(1-|z|^2)^2}\delta_{\alpha\bar\beta}$. 
Then Assumption \ref{mainassumption} is valid for $\mu>1$ and $\mu':=\mu-1$.
\end{example}

\begin{proof}Since $\gamma$ is rotational symmetric, Lemma \ref{lines} states that geodesics through $0$ are lines, hence the
distance function is given by
$$d(0,z)=\int_0^{|z|}\frac {2}{1-t^2}dt=2\mbox{artanh}(|z|).$$

For a rotational symmetric function $f(r)$ on ${\bb R}^4$ the ordinary Laplacian is given by
$$\Delta f=\left(\frac{\partial^2}{\partial r^2}+\frac 3r\frac{\partial}{\partial r}\right)f.$$
So we compute
\begin{eqnarray*}\Delta(A+2\mbox{artanh}(r))^{-\mu'}&=&\frac{4\mu'(\mu'+1)}{(1-r^2)^2}(A+2\mbox{artanh}(r))^{-(\mu'+2)}-
\frac{4\mu' r}{(1-r^2)^2}(A+2\mbox{artanh}(r))^{-(\mu'+1)}\\
& &-\frac{6\mu'}{r(1-r^2)}(A+2\mbox{artanh}(r))^{-(\mu'+1)}
\end{eqnarray*}
and hence
\begin{eqnarray*}-\tilde\Delta(A+2\mbox{artanh}(r))^{-\mu'}&=&-\frac 14(1-r^2)^2
\Delta(A+2\mbox{artanh}(r))^{-\mu'}\\
&=&-\mu'(\mu'+1)(A+2\mbox{artanh}(r))^{-(\mu'+2)}+
\mu' r(A+2\mbox{artanh}(r))^{-(\mu'+1)}\\
& &+\frac{3\mu'(1-r^2)}{2r}(A+2\mbox{artanh}(r))^{-(\mu'+1)}
\end{eqnarray*}

Elementary calculations show that the coefficient of $(A+2\mbox{artanh}(r))^{-(\mu'+1)}$ is strictly decreasing,
$$\frac{\partial}{\partial r}\left(\mu' r+\frac{3\mu'(1-r^2)}{2r}\right)<0$$
and hence
$$\mu' r+\frac{3\mu'(1-r^2)}{2r}>\mu'.$$
If we now choose $A>\mu'+1$, then
$$-\tilde\Delta(A+2\mbox{artanh}(r))^{-\mu'}>C(A+2\mbox{artanh}(r))^{-\mu'-1},$$
with $C:=\frac{A-\mu'-1}{A}$.

With arguments as above this implies that Assumption (1) is satisfied for any $\mu>1$ and $\mu':=\mu-1$.  
\end{proof}

\vspace{4mm}
\noindent
Now we turn our attention to the norm of the tension field. We will see that this is the crucial point.

\begin{proposition}\Punkt If $M$ is a complex manifold with Hermitian metrics $\gamma$ and $\tilde\gamma$, 
and if $id:(M,\gamma)\seq (M,\tilde\gamma)$
denotes the identity map, then we define
$$A^{\varepsilon}:=\frac{1}{2}
      \tilde{\gamma}^{\varepsilon\bar{\delta}}\gamma^{\alpha\bar{\beta}}
\left(\tilde{\gamma}_{\alpha\bar{\delta},\bar\beta}
   -\tilde{\gamma}_{\alpha\bar\beta,\bar\delta}
\right).$$
With this vector field given we obtain
$$\|\sigma(id)\|^2=\gamma_{\varepsilon\bar\phi}{A^\varepsilon}\overline{A^\phi}.$$
\end{proposition}

\vspace{4mm}
\noindent
This formula simplifies in the conformal case:

\begin{proposition}\Punkt If $\gamma=f\delta_{\alpha\bar\beta},\tilde\gamma=\tilde f\delta_{\alpha\bar\beta}$, 
with smooth real valued positive functions $f, \tilde{f}$, then
$\|\sigma(id)\|=\frac{m-1}{2 f}\left|\nabla\sqrt{\tilde{f}} \right|$. 
In particular,
if $\gamma=\tilde\gamma=f\delta_{\alpha\bar\beta}$, then 
$\|\sigma(id)\|=\frac{(m-1)}{2}\left|\nabla\frac{1}{\sqrt f}\right|$.
\end{proposition}

\begin{example}\Punkt In particular, if $M=B=B_1 (0)\subset \mathbb{C}^m$ and $\gamma=\tilde\gamma
=\frac{4}{(1-|z|^2)^2}\delta_{\alpha\bar\beta}$ 
is the Poincar\'e metric, then
$$\|\sigma(id)\|=(m-1)|z|.$$
\end{example}

\vspace{4mm}\noindent
So $\|\sigma(id)\|\not\in C^0_\mu(M)$ for the Poincar\'e case, we do not even have decay to zero. 
We set this result in
a more general framework now.

Let $\Omega\subset{\bb C}^m$
be a rotational symmetric domain, equipped with a rotational symmetric metric $\gamma$, which obeys 
the polar block form condition of
Lemma \ref{block}, if
$0\notin\Omega$. Let 
$S\subset\Omega$ be a sphere centered in $0$ with radius $r_0$. Then, for $r\ge 0$ choose $x\in\Omega$ with $|x|=r$ and define
$$D(r):=\mbox{dist}(x,S).$$
Obviously, $D(r)$ is well-defined and for any fixed $s\in S$ the function $D(r)$ obeys 
$$D(r)\le d(x,s)\le D(r)+C,$$
where $C$ depends only on $b_{ij}(r_0)$. Hence statements about growth of $d(x,s)$ are equivalent to those about $D(r)$ and
independent of the choices of $S$ and $s$. A calculation similar to that in the proof of Lemma \ref{block}
shows that $\mbox{dist}(x,S)$ is realized by the segment of the line containing $x$ and $0$.
If $0\in\Omega$ we set $r_0:=0$, i.e. $D(r)=d(r,0)$.

\begin{proposition}\Punkt Let $\Omega\subset{\bb C}^m$ be equipped with a rotational symmetric, complete metric $\gamma$ 
conformal to the euclidean such that $\|\sigma(id)\|\in C^0_\lambda(\Omega)$ for $\lambda>1$. Then $D(r)$ has linear growth.
In particular, $\Omega={\bb C}^m$.
\end{proposition}

\begin{proof}We denote $\gamma=f\delta_{\alpha\bar\beta}$ and 
abbreviate $h:=\frac{1}{\sqrt f}$. We calculate
$$\|\sigma(id)\|=(m-1)|h'|=(m-1)\left|\left(\frac 1{D'}\right)'\right|
=(m-1)\left|\frac{D''}{(D')^2}\right|<CD^{-\lambda}$$
for $D\gg 0$. Since $D$ is strictly increasing for $r>r_0$, this implies for $D\gg 0$
$$\left|\frac{D''}{D'}\right|<C_1 D^{-\lambda}D'.$$
Integration yields
$$\mbox{Var }(\ln D',[r_0+\varepsilon,r])<C_3-C_2D^{1-\lambda}.$$
If $r>r_0+\varepsilon$ increases, $\mbox{Var }(\ln D',[r_0+\varepsilon,r])$ is increasing, as well as
$D(r)$, hence $$\mbox{Var }(\ln D',\{r>r_0+\varepsilon\})\le C_3.$$
This implies
$$0<C_4<D'(r)<C_5$$
for all $r\ge r_0$ and hence
$$C_4r+C_6<D(r)<C_5r+C_7.$$
Very similar arguments apply for $r<r_0$, if $0\notin\Omega$. Hence $0\in\Omega$ and $\Omega={\bb C}^m$.
\end{proof}

\begin{proposition}\Punkt Any rotational symmetric Hermitian metric $\gamma$ on ${\bb C}^m$, 
which has nonpositive sectional curvature is
either euclidean or $D(r):=d(0,r)$ has superlinear growth.\end{proposition}

\begin{proof}First, we reduce to the case $m=1$: If $M$ has nonpositive sectional curvature and $E$
is a complex plane through $0$, then $M\cap E$ has also nonpositive sectional curvature. On the other hand,
if $\gamma|_{M\cap E}$ is euclidean or has superlinear growth for some plane $E$ containing $0$, we conclude by the
rotational symmetry that this holds also for $\gamma$. So we assume now $m=1$.

With notation as above, we compute for $\gamma=f(r)dz\otimes d\bar z=\phi(r^2)dz\otimes d\bar z$
and $s:=r^2$
$$\phi(s)R_{1212}=-\left[2s\left(\phi''(s)\phi(s)-(\phi'(s))^2\right)+2\phi'(s)\phi(s)\right].$$
This implies
$$2\phi^2(s(\ln\phi)''+(\ln\phi)')\ge 0.$$
We abbreviate $(\ln\phi)'=:\psi$. Since $R=-\frac{\Delta \ln f}{2f}$ the
maximum principle implies that $f(r)$ is increasing and hence $f'(r)\ge 0$ and
also $\phi'(s)\ge 0$. So we conclude $$\psi\ge 0.$$ 

We claim that $\psi(s)>0$ for $s>0$ unless $\psi\equiv 0$. To prove this we assume that there are $0<s_1<s_2$ such that
$\psi(s_1)>0$ and $\psi(s_2)=0$. Then
$$0\le\int_{s_1}^{s_2}(s\psi'(s)+\psi(s))ds=s_2\psi(s_2)-s_1\psi(s_1)=-s_1\psi(s_1)<0,$$
what is a contradiction. Hence $\psi(s)>0$ for all $s>0$ or $\psi\equiv 0$. The last case is the euclidean case.

So we assume $\psi(s)>0$ for $s>0$. Let us fix some $s_0>0$. 
 
Then we compute for all $s>s_0$
\begin{eqnarray*}s\psi'+\psi\ge 0  & \iff & \frac{\psi'}{\psi}\ge-\frac 1s\\
&\imp &\ln\psi\ge C_1-\ln s\\
&\imp &(\ln\phi)'\ge \frac{C_2}{s}\mbox{ with }C_2>0\\
&\imp &\phi(s)\ge C_3s^{C_2}\mbox{ with }C_3>0\\
&\imp &f(r)\ge C_3r^{2C_2}\\
&\imp &D(r)\ge C_4+C_5r^{1+C_2}\mbox{ with }C_2,C_5>0
\end{eqnarray*}
In the computations we always integrate from $s_0$ to $s$. 
\end{proof}

\begin{cor}\Punkt\label{eucl}If $M\subset{\bb C}^m$ allows for a rotational symmetric complete metric conformal to the euclidean
with nonpositive sectional curvature
and $\|\sigma(id)\|\in C^0_\lambda(M)$ for some $\lambda>1$, then $M={\bb C}^m$ and the metric is the
euclidean metric multiplied with a constant. In particular, $\sigma(id)\equiv 0$.
\end{cor}

\vspace{4mm}
\noindent
These results illustrate that there is no obvious example, where the identity
map $id$ may serve as initial map $h$.

In order to construct nontrivial Hermitian-harmonic maps, one might look for
manifolds with \emph{two} infinite ends as in Example~\ref{zweienden} above.
However, if one wants to choose $M=N$ in this case, one has to observe the
following obstruction:

\begin{rem}\Punkt
$N={\bb R}^n\setminus\{0\}$ does not admit a complete, nonpositively curved metric for $n\ge 3$.
If we would have a nonpositively curved metric, the Cartan-Hada\-mard-theorem would imply that
${\bb R}^n$ is the universal cover of $N$. Since $N$ is simply connected for $n\ge 3$, $N$ would have to be isomorphic
to ${\bb R}^n$. But since $\pi_{n-1}({\bb R}^n)=0$ and $\pi_{n-1}(N)={\bb Z}$, this is not the case.
\end{rem}

\section{The corresponding parabolic system}\label{heatequation}

\subsection{Existence and convergence results}\label{existenceconvergence}

Originally in the fundamental work of Jost and Yau \cite{JostYau},
as in many contributions to the harmonic map system, existence results
were obtained via the seeming detour of the corresponding heat system
\begin{equation}\label{heatsystem}
\left\{\begin{array}{ll}
\displaystyle \frac{\partial u}{\partial t} = \sigma (u)
                        &\mbox{\ on\ } (0,\infty) \times M,\\[2mm]
\displaystyle u(0) = h,\\[2mm]
\displaystyle u \sim h & \mbox{\ at infinity, homotopic to each other}.
\end{array}\right.
\end{equation}
The initial map is chosen as in Theorem~\ref{existenceresult}, and here,
the notation ``initial map'' as well as the homotopy between $h$ and $u$
become more transparent.

Similarly as in Theorem~\ref{existenceresult} we get existence of a global solution
to (\ref{heatsystem}) and also convergence to {\it some smooth map}
$u: M \to N$ for a sequence $t_k \to \infty$
by means of an exhaustion procedure.

\newpage
\begin{theorem}[Global Existence]\Punkt\label{globalexistence}
Assume that $M$ is a noncompact complete Hermitian manifold such that for the
holomorphic Laplace operator $- \tilde{\Delta} $ on $M$, the
Assumption~\ref{mainassumption} is satisfied with positive numbers
 $\mu,\mu'>0$. Further let $N$ be a
complete Riemannian manifold with nonpositive sectional curvature
and $h:M \to N$ a smooth map with $\| \sigma (h) \|  \in  C^0_{\mu} (M)$.

Then there exists a global smooth solution
$u: [0, \infty) \times M \to N$ to (\ref{heatsystem})
such that $u(t,\, . \, )$ is for every $t \ge 0$ homotopic
to $h$. For the homotopy distance $\rho(t,\, . \, )$ between $h(\, . \, )$
and $u(t,\, . \, )$, we have that $\rho(t) \in  C^0_{\mu'} (M)$
uniformly in $t$.

Moreover there exists a sequence $t_k \to \infty$ such that
$u(t_k,\, . \, )$ converges to a smooth map $u:M \to N$,
being homotopic to $h$ and converging to $h$ at ``infinity''.
\end{theorem}

\vspace{4mm}\noindent
{\it Proof.}  As above, let $\Omega_k$ be  a compact exhaustion of $M$.
According to \cite[Proof of Theorem 6]{JostYau} , there exist
smooth solutions 
$u_k: [0,\infty) \times \Omega_k \to N$ of the
initial boundary value problems
\begin{equation}\label{compactIBVProblems}
\left\{ \begin{array}{ll}
\displaystyle (u_k)_t -\sigma(u_k) = 0  &\mbox{\ in\ }
                       [0,\infty) \times \Omega_k,\\[2mm]
\displaystyle u_k (t,x) = h(x)           &\mbox{\ for\ }
           (t,x) \in     [0,\infty) \times    \partial\Omega_k,\\[2mm]
\displaystyle u_k (0,\, . \, ) = h           &\mbox{\ on\ }
                     \Omega_k,\\[2mm]
\displaystyle u_k \mbox{\ homotopic to\ } h, &\mbox{\ with respect to\ }
 \partial\Omega_k.
\end{array} \right.
\end{equation}
We first need to show that $\left\| \frac{\partial u_k}{\partial t} (t,x) 
\right\|$
are uniformly bounded. For this purpose we note that according to
formula \cite[(6.2)]{LeiNi},
$\left\| \frac{\partial u}{\partial t} (t,x) \right\|$
satisfies an initial boundary value problem
for the following differential inequality
\begin{equation}\label{decay}
\left\{ \begin{array}{ll}
\displaystyle \left( \frac{\partial}{\partial t}
   - \frac{1}{4} \tilde{\Delta}  \right) 
  \left\| \frac{\partial u_k}{\partial t} (t,x) \right\| \le 0  &\mbox{\ in\ }
                       [0,\infty) \times \Omega_k,\\[4mm]
\displaystyle  \left\| \frac{\partial u_k}{\partial t} (t,x) \right\| = 0
  &\mbox{\ for\ } (t,x) \in     [0,\infty) \times    \partial\Omega_k\\[4mm]
\displaystyle  \left\| \frac{\partial u_k}{\partial t} (0,x) \right\|
            = \left\| \sigma (h) (x) \right\|        &\mbox{\ for\ } 
                    x \in  \Omega_k.
\end{array} \right.
\end{equation}
By means of the parabolic maximum principle, we conclude that
for any $k$ and all $(t,x) \in [0,\infty) \times \Omega_k$ we have
\begin{equation}
\left\| \frac{\partial u_k}{\partial t} (t,x) \right\|
 \le \max_{x \in M} \| \sigma (h(x)) \|.
\end{equation}
Next we need a $L^{\infty}$ bound for $\rho_k$ on any $[0,T] \times
\Omega_k$, where again $\rho_k(t,x)$ denotes the homotopy distance between
$u_k (t,x)$ and $h(x)$. For this purpose we again introduce
a nonnegative barrier function
$V \in  C^0_{\mu'} (M)$, such that
\begin{equation}
- \tilde{\Delta}  V = 4 \| \sigma (h) \| \mbox{\ in\ } M.
\end{equation}
Furthermore we note that the crucial estimate
(\ref{Corollary3.5}) generalizes to smooth time dependent maps
$u,v: [0,\infty) \times M \to N$ as follows:
\begin{equation}
\left( \frac{\partial}{\partial t}- \frac{1}{4} \tilde{\Delta}\right) 
\rho(u,v) \le
       \left\| \frac{\partial u}{\partial t} - \sigma (u)\right\| 
     +  \left\| \frac{\partial v}{\partial t} - \sigma (v) \right\|.
\end{equation}
We conclude
\begin{equation}
\left\{ \begin{array}{ll}
\displaystyle \left( \frac{\partial}{\partial t} 
   - \frac{1}{4} \tilde{\Delta}  \right) \rho_k
  \le \| \sigma (h) \|
  \le  \left( \frac{\partial}{\partial t}
   - \frac{1}{4} \tilde{\Delta}  \right) V
                 &\mbox{\ in\ }
                       [0,\infty) \times \Omega_k,\\[4mm]
\displaystyle  \rho_k (t,x) = 0    \le V (x)
  &\mbox{\ for\ } (t,x) \in     [0,\infty) \times    \partial\Omega_k,\\[4mm]
\displaystyle  \rho_k (0,x) = 0 \le V(x)     &\mbox{\ for\ }
                    x \in  \Omega_k.
\end{array} \right.
\end{equation}
By means of the parabolic maximum principle, we get
for all $(t,x) \in [0,\infty) \times \Omega_k$
\begin{equation}
\rho_k (t,x) \le V (x),
\end{equation}
i.e. the desired uniform $L^{\infty}$ bound for the $\rho_k$
and hence for the $\rho_{kj}=\rho(u_k, u_j)$.
Analyzing and correcting carefully the argument 
given in \cite[pp. 351--352]{LeiNi},
one obtains for any $T> 0$ and any relative compact $\Omega \subset M$
a uniform (in $k$ and $T$) bound for
$$
\int_T^{T+2} \int_\Omega e(u_k) .
$$
{From} this we want to deduce a local  maximum bound for
$e(u_k)$, which by means of standard linear parabolic theory 
will allow to pass to the limit and to obtain a
{\it global} smooth solution
to (\ref{heatsystem}). For this purpose we may assume that
$\Omega$ is contained in one single coordinate chart, and we
take from \cite[(6.5)]{LeiNi} that $e(u_k)$ satisfies
a differential inequality of the form
\begin{equation}
\left( \frac{\partial}{\partial t} -
\frac{1}{4}  \tilde{\Delta} + C(\Omega) \right)
   e(u_k) \le 0.
\end{equation}
The constant $C(\Omega)$ may be suitably chosen independently
of $k$. Here we have to apply the local maximum principle
for parabolic operators not in divergence form, which can be
adapted from \cite[Theorem 7.21]{Lieberman}.
(For an extensive discussion we refer to 
Proposition~\ref{parmax} in the appendix.)
This gives a bound
for $\max_{[T+1,T+2] \times B} e(u_k)$ for  sufficiently small balls
contained in $\Omega$, which depends on $\max e(h)$, $B$, $\Omega$
and the $L^1$-bound on $e(u_k)$, but not on $T$.
Hence we have found a maximum bound for the gradient, which is
local in space, but {\it global in time}.

Homotopy between $u$ and $h$ is shown as in the proof
of \ref{existenceresult}. Moreover we note that the homotopy distance
is also bounded by $V \in  C^0_{\mu'} (M)$:
\begin{equation}
\rho(t,x) := \rho (u(t,x), h(x)) \le V(x).
\end{equation}
The stated convergence now follows from the (uniform in time, local
in space) boundedness of $e(u_k)$ and the global boundedness
of $\left\| \frac{\partial u}{\partial t} (t,x) \right\|$
by standard linear parabolic theory.
\ \hfill $\square$

\vspace{4mm}\noindent
Next, we want to prove that this global solution converges to
a (stationary)  Hermitian-harmonic map $u:M \to N$. Here it seems
that we need something stronger than Assumption~\ref{mainassumption}.
As additional hypothesis we formulate:

\begin{assumption}[Decay properties]\Punkt\label{secondassumption}
\ \ \ \\
We assume that there exists a positive number $\mu>0$ such that
for every $\varphi \in C^0_{\mu} (M)$, we have  decay
of $\max v(t,\, . \, ) $ towards $0$ for every bounded solution $v$
of the initial value problem for the heat equation with the holomorphic
Laplace operator and $\varphi$ as initial datum. Moreover we
assume that the solution of the initial value problem is unique
in the class of all uniformly bounded functions on $[0,T] \times M$.
\end{assumption}

\vspace{4mm}
\noindent
The formulation of Assumption \ref{secondassumption} suggests the use of a comparison function, what is a great difficulty
for arbitrary metrics. Hence it would be more adequate to find a spectral reformulation. This is aimed at by the following
Lemma. We denote
$S_{\phi_1,\phi_2}:=\{r\exp(i\phi)|\phi_1\le\phi\le\phi_2\}$ and $C^0_b(M):=C^0_{\mu=0}(M)$ for the bounded continuous functions
in order to avoid confusion with compactly supported functions.

\begin{lemma}\Punkt
Assume that $(-\tilde\Delta+\lambda)f=0$ has a unique bounded solution for every initial datum $\varphi\in C^0_\mu(M)$ and
moreover
$$\|(-\tilde\Delta+\lambda)^{-1}\|_{B(C^0_\mu(M),C^0_b(M))}\le C$$
uniformly for all $\lambda\in S_{\phi_1,\phi_2}$ for certain $\pi/2<\phi_1<\pi, -\pi/2>\phi_2>-\pi$. 
Then Assumption \ref{secondassumption} holds.
\end{lemma}

\begin{proof}
The key ingredient is the keyhole integral. Let $\Gamma=\Gamma_1+\Gamma_2
+\Gamma_3$ be a path in ${\bb C}$ such that
\begin{eqnarray*}\Gamma_1&:=&\{r\exp(i\phi_1)|r\in\left[1;\infty\right)\}\\
\Gamma_2&:=&\{\exp(i\phi)|\phi\in[\phi_1;\phi_2]\}\\
\Gamma_3&:=&\{r\exp(i\phi_2)|r\in\left[1;\infty\right)\}.\end{eqnarray*}
Then we can define a semigroup (cf. e.g. \cite[Part II]{Friedman})
$$\exp(\tilde\Delta t):=\frac 1{2\pi i}\int_\Gamma\exp(\lambda t)(-\tilde\Delta+\lambda)^{-1}d\lambda$$
as operator in $B(C^0_\mu(M),C^0_0(M))$, since we assumed the uniform boundedness of $(-\tilde\Delta+\lambda)^{-1}$.
It has the property
$$\frac{d}{dt}(\exp(\tilde\Delta t) u)=\tilde\Delta\exp(\tilde\Delta t)u.$$
Cauchy's integral formula implies that integration over $\Gamma$ yields the same as integration over $\Gamma/t$ for
$t>0$. Hence
\begin{eqnarray*}
\|\exp(\tilde\Delta t)\|&=&\frac 1{2\pi}\left\|\int_{\Gamma/t} \exp(\lambda t)
\left(-\tilde\Delta+\lambda\right)^{-1}d\lambda\right\|\\
&=&\frac 1{2\pi t}\left\|\int_\Gamma\exp(\lambda)
\left(-\tilde\Delta+\frac{\lambda}{t}\right)^{-1}d\lambda\right\|\\
&\le&\frac 1{2\pi t}\int_\Gamma   | \exp(\lambda) | \left\|
\left(-\tilde\Delta+\frac{\lambda}{t}\right)^{-1}\right\|\, |d\lambda|   \\
&\le&\frac C{2\pi t}\int_\Gamma | \exp(\lambda) | \, \left| d\lambda\right|\\
&=&\frac{C'}{t},
\end{eqnarray*}
since $\cos(\phi_1)<0$ and $\cos(\phi_2)<0$.

Hence the unique bounded solution $v:=\exp(\tilde\Delta t)\varphi$ satisfies
$$\max_{x\in M}v(t,x)\le\frac{C'\|\varphi\|}{t}.$$
So Assumption \ref{secondassumption} is valid.
\end{proof}

\vspace{4mm}\noindent
The meaning  and relevance of Assumption \ref{secondassumption} will be extensively
discussed in the examples 
in subsection \ref{parabolicexamples} below. With help of this assumption, we may now state:

\begin{theorem}[Convergence to a Hermitian-harmonic map]\Punkt\newline
Let the assumptions of Theorem~\ref{globalexistence} be satisfied
as well as Assumption~\ref{secondassumption}
with the same $\mu$ as in Assumption~\ref{mainassumption}. 
Then, for the solution
$u(t, \, . \, ) $ of the time dependent Hermitian-harmonic map system
(\ref{heatsystem}), there exists a sequence $t_k \to \infty$ such that
$u(t_k, \, . \, ) $ converges to a Hermitian harmonic map $u:M \to N$.
\end{theorem}

\vspace{4mm}\noindent
{\it Proof.} It remains to show a decay result for $\max_M \left\| 
\frac{\partial u}{\partial t} \right\|$. The latter is achieved
by means of the differential inequality
\begin{equation}
\left\{ \begin{array}{ll}
\displaystyle \left( \frac{\partial}{\partial t} 
   - \frac{1}{4} \tilde{\Delta}  \right) 
  \left\| \frac{\partial u}{\partial t} (t,x) \right\| \le 0  &\mbox{\ in\ }
                       [0,\infty) \times M,\\[4mm]
\displaystyle  \left\| \frac{\partial u}{\partial t} (0,x) \right\|     
            = \left\| \sigma (h) (x) \right\|        &\mbox{\ for\ }
                    x \in  M.
\end{array} \right.
\end{equation}
Since $\left\|  \frac{\partial u}{\partial t} \right\|$ is uniformly
bounded, Assumption~\ref{secondassumption} gives decay to $0$,
as $t \to \infty$.\ \hfill $\square$

\subsection{Examples}\label{parabolicexamples}

In the remainder, we show that Assumption~\ref{secondassumption} is
likewise satisfied in all the examples treated above in
Subsection~\ref{examples}.

\begin{example}\Punkt Let $M={\bb C}^m$ and $\gamma=\delta_{\alpha\bar\beta}$ be the euclidean metric. 
Then Assumption \ref{secondassumption} holds true. 
\end{example}

\begin{proof}
Assume
$\varphi\in C^0_\mu(M)$ with $\mu\le 2m$, for simplicity we specialize to $|\varphi(y)|<(1+|y|)^{-2\alpha m-\varepsilon}$ with 
$\alpha\in(0;1]$ and $\varepsilon>0$ such that $2\alpha m+\varepsilon =\mu$.
Since the fundamental solution of the heat equation is
$$\gamma(t,x)=C_0t^{-m}\exp\left(-\frac{|x|^2}{4t}\right),$$
the solution $v(t,x)$ of the initial data problem with $\varphi$ as initial datum is given by
\begin{eqnarray*}|v(t,x)|&=&C_0 
 \left| t^{-m}\int\exp\left(-\frac{|x-y|^2}{4t}\right)\varphi(y)dy \right|\\ 
&\le&
C_0 t^{-m}\left(\int\exp\left(-\frac{|x-y|^2}{4t}
\right)^{\frac{1}{1-\alpha}}dy\right)^{1-\alpha}
\left(\int|\varphi(y)|^{\frac 1\alpha}dy\right)^\alpha\\ 
&\le&
C_1 (\varphi)
t^{-m}\left(\int\exp\left(-\frac{|x-y|^2}{4t}\right)^{\frac{1}{1-\alpha}}dy
\right)^{1-\alpha},\\ 
& &\mbox{since }|\varphi|^{\frac 1\alpha}\mbox{ is
integrable,}\\ 
&\le& C_1  (\varphi)   
t^{-m}\left(\int\exp\left(-\frac{|x-y|^2}{4t}\right)dy\right)^{1-\alpha}\\ & =
& C_2(\varphi)
t^{-m}\left(t^m\int\exp(-|z|^2)dz\right)^{1-\alpha}\\ 
& = & C_3(\varphi)  t^{-\alpha m}.
\end{eqnarray*}

If $\mu>2m$, then $\varphi$ is integrable and hence
$$|v(t,x)|\le C_4t^{-m}.$$

This proves the validity of Assumption~\ref{secondassumption} for the euclidean case for every $\mu>0$.
\end{proof}

\vspace{4mm}
\noindent
To be able to treat the case of Example \ref{basicex}, we have to formulate a maximum principle for the corresponding
Laplace operator.

\begin{lemma}\Punkt\label{confmax} Let $M={\bb C}^m$ and $\gamma=(1+r^2)^{-1}\delta_{\alpha\bar\beta}$. If
$u:M\times [0;T]\seq{\bb R}$ is bounded and satisfies $(-\tilde\Delta+\frac{\partial}{\partial t})u\ge 0,u(x,0)\ge 0$, then $u\ge 0$.
\end{lemma}

\begin{proof}We imitate the proof of the euclidean case like given in \cite[V,4,Thm4.1]{Di}. We take the function
$$v(x,t):=(1+|x|^2)\exp(4mt)\ge 1+|x|^2,$$
which also satisfies the heat equation, i.e. $(-\tilde\Delta+\frac{\partial}{\partial t})v=0$.
For the function
$$w_\varepsilon:=u+\varepsilon v$$
we hence get
$$(-\tilde\Delta+\frac{\partial}{\partial t})w_\varepsilon\ge 0, w_\varepsilon(x,0)\ge 0$$
and $w_\varepsilon(x,t)\ge 0$ outside a compact set $K_\varepsilon\times [0;T]\subset M\times [0;T]$.
Now using the parabolic maximum principle for the bounded domain $K_\varepsilon\times [0;T]$ we see that
$w_\varepsilon\ge 0$ everywhere. Hence $u(x,t)=\lim_{\varepsilon\seq 0}w_\varepsilon(x,t)\ge 0$.
\end{proof}

\vspace{4mm}
\noindent
Now we are able to continue Example \ref{basicex}.

\begin{example}\Punkt
Let us consider the conformal metric $\gamma_{\alpha\bar\beta}=(1+r^2)^{-1}\delta_{\alpha\bar\beta}$ on $M={\bb C}^m$, $m\ge 2$. 
Then again Assumption \ref{secondassumption} holds true.
\end{example}

\begin{proof}
Recall that the geodesic length $d(0,x)\sim \ln(1+r^2)$.
The decay of the
solution of the corresponding heat equation is proven by comparison to a test function. Let us consider
$$w(t,x):=(A+\ln(1+r^2)+t)^{-\mu}.$$
Now we verify
$$
\left(-\tilde\Delta+\frac{\partial}{\partial t} \right)w\ge 0,
$$
if $A$ is chosen big enough, given $\mu$ and $m$.

We compute
\begin{eqnarray*}\left(-\tilde\Delta+\frac{\partial}{\partial t}\right)w 
&=&\mu\Big(((4m-5)A-4\mu-4){r}^{2}+(4m-5){r}^{2}\ln (1+{r}^{2})+(4m-5)
{r}^{2}t\\
& &\ \ \ +(4m-1)A+(4m-1)\ln (1+{r}^{2})+(4m-1)t\Big)\\
&&\cdot \left(1+\ln
(1+{r}^{2})+t\right )^{-\mu-2} \left (1+{r}^{2}\right )^{-1}
\end{eqnarray*}

Obviously, $\left(-\tilde\Delta+\frac{\partial}{\partial t}\right)w>0$, if
$$(4m-5)A-4\mu-4\ge 0.$$
So, given $\mu$ and $m\ge 2$ we choose $A\ge 1$ such that this inequality is satisfied. 

Since for an exact bounded solution with $|v(0,x)|\in C^0_\mu(M)$ we have $|v(0,x)|\le Cw(0,x)$, 
the parabolic maximum principle
stated in Lemma \ref{confmax}
proves
$$|v(t,x)|\le Cw(t,x)\le Ct^{-\mu}.$$

The uniqueness in the class of bounded solutions immediately follows by Lemma \ref{confmax} and the observation that
$u$ is a solution of the heat equation with zero initial data if and only if $-u$ is.

Hence Assumption~\ref{secondassumption} is valid for all $\mu>0$.
\end{proof}

\vspace{4mm}
\noindent
Finally, we come back to Example \ref{zweienden}. Again we first have to prove a maximum principle.

\begin{lemma}\Punkt\label{max4}Let $M=({\bb S}^1)^{2m-1}\times(-1,1)$ with 
$\gamma=\delta^2(1-|s|)^{-2\delta-2}\delta_{\alpha\bar\beta}$ 
for $1/2<|s|<1$ like in Example \ref{zweienden}. If $u:M\times[0;T]\seq{\bb R}$
is bounded and satisfies $(-\tilde\Delta+\frac{\partial}{\partial t})u\ge 0$ and
$u(x,0)\ge 0$, then $u\ge 0$.
\end{lemma}

\begin{proof}Like in the proof of Lemma \ref{confmax} it is sufficient to construct a supersolution $v(s,t)$ such
that $\inf_{t\in[0;T]}v(s,t)\seq\infty$ for $s\seq\pm 1$. The choice
$$\tilde v(s,t):=(1-\log(1-|s|))\exp(\frac{t}{\delta^2})$$
works for $|s|>\frac 12$. Since for $|s|\le \frac 12$ there exists $C=C(T)>0$
such that $$\left(-\tilde\Delta+\frac{\partial}{\partial t}\right)\tilde v\ge
-C,$$ the function $v(s,t):=\tilde v(s,t)+Ct$ is a supersolution satisfying
the required conditions. 
\end{proof}

\begin{example}\Punkt Let $M$ be like in Example \ref{zweienden}. Then
Assumption \ref{secondassumption} is valid.
\end{example}

\begin{proof}
Let us define $(a)_i:=a(a+1)\cdot ...\cdot(a+i-1), (a)_0:=1$ for real $a$ and integer $i$. 
Denote the Kummer function
$$F(a,b,z):=\subs F11(a,b,z)=\sum_{i=0}^\infty\frac{(a)_i}{(b)_ii!}z^i,$$
which is convergent for all $z$, if $b$ is not a negative integer.
We will make use of the following properties:

\begin{eqnarray}0&=&z\frac{\partial^2}{\partial z^2}F(a,b,z)+(b-z)\frac{\partial}{\partial z}F(a,b,z)-aF(a,b,z)
\label{diffequ}\\
F(a,b,z)&=&\frac{\Gamma(b)}{\Gamma(a)}\exp(z)z^{a-b}(1+O(|z|^{-1}))\mbox{ if }\Re(z)>0\label{limit}\\
F(a,b,z)&=&\exp(z)F(b-a,b,-z)\label{kummer}\\
aF(a+1,b,z)&=&aF(a,b,z)+z\frac{\partial}{\partial z}F(a,b,z)\label{diff}
\end{eqnarray}

All these properties can be found in \cite{abr} as $13.1.1, 13.1.4, 13.1.27, 13.4.10$.

As a comparison function we choose
$$
w(s,t):=t^cF \left(-c,1+\frac{1}{2\delta},
   -\frac{1}{4}(1-|s|)^{-2\delta}t^{-1}\right),
$$
where we choose $\max(-\frac{1}{2\delta},-\frac{\mu}{2})<c<0$.
By the first Kummer transformation (\ref{kummer}) this becomes
$$
w(s,t)=t^c\exp(-\frac
1{4}(1-|s|)^{-2\delta}t^{-1})F\left(1+\frac{1}{2\delta}+c,
1+\frac{1}{2\delta},\frac{1}{4}(1-|s|)^{-2\delta}t^{-1}\right),
$$ 
and now we
see easily that $w(s,t)>0$. By (\ref{limit}) for fixed $s$ and $t\seq 0$ the
function $w(s,t)$ is continuously extendable to $t=0$ and $$
w(s,0)=C(1-|s|)^{-2c\delta}>C(1-|s|)^{\delta\mu}=C(1+\tilde d(s))^{-\mu}.
$$
By a simple calculation using (\ref{diffequ}) we can see that $w(s,t)$ is an exact solution to
\begin{equation}
\left(-\delta^{-2}(1-|s|)^{2\delta+2}\Delta+
\frac{\partial}{\partial t}\right)w(s,t)=0\label{approxheat}
\end{equation}
on $M$ outside $s=0$. 
Since this is not yet the original equation and $w$ is singular in $s=0$ we have to 
do some more calculations. First we note that
by (\ref{diff}) and (\ref{kummer})
\begin{eqnarray*}
\frac{\partial}{\partial t}w(s,t)&=&ct^{c-1}
F \left(-c+1,1+\frac{1}{2\delta},
    -\frac{1}{4}(1-|s|)^{-2\delta}t^{-1}\right)\\
&=&ct^{c-1}\exp(-\frac 1{4}(1-|s|)^{-2\delta}t^{-1})
F \left(\frac{1}{2\delta}+c,1+\frac{1}{2\delta},
\frac{1}{4}(1-|s|)^{-2\delta}t^{-1}\right)\\
&<&0,
\end{eqnarray*}
since $c<0, \frac{1}{2\delta}+c>0$ and $F(a,b,z)$ is positive, if $a,b,z>0$. 
{From} this and the fact that $w(s,t)$ is a solution of (\ref{approxheat}) we deduce that
$$-\Delta w(s,t)>0$$
for all $0\not=s,t> 0$. If we determine $C>0$ such that
$Ca(s)<\delta^2(1-|s|)^{-2\delta-2}$ for $|s|\le \frac 12$  then $\tilde
w(s,t):=w(s,Ct)$ satisfies 
$$
\left(-\tilde\Delta+\frac{\partial}{\partial t}\right)\tilde
w(s,t)\ge 0
$$ 
for $s\not=0$. So, for the sake of simplicity let us assume that
$C=1$.

Now we are considering $s>0$. Then, again by (\ref{diff}) and (\ref{kummer})
\begin{eqnarray*}\frac{\partial}{\partial s}w(s,t)&=&-2\delta ct^c(1-s)^{-1}
\left(F  \left(-c+1,1+\frac{1}{2\delta},
-\frac{1}{4}(1-|s|)^{-2\delta}t^{-1}\right)\right.\\
& &\left.
-F\left(-c,1+\frac{1}{2\delta},
      -\frac{1}{4}(1-|s|)^{-2\delta}t^{-1}\right)\right)\\
&=&2\delta ct^c(1-s)^{-1}\exp(-\frac 1{4}(1-s)^{-2\delta}t^{-1}) \left(F
\left(1+\frac{1}{2\delta}+c,1+\frac{1}{2\delta},    
\frac{1}{4}(1-s)^{-2\delta}t^{-1}\right)\right.\\  
& &\left. -F
\left(\frac{1}{2\delta}+c,1+\frac{1}{2\delta},     
\frac{1}{4}(1-s)^{-2\delta}t^{-1}\right)\right)\\ &<& 0\end{eqnarray*}
again by $c<0$ and the property that $F(a,b,z)>F(a',b,z)$ if $a>a'>0, b,z>0$.

Finally, with this in mind we are going to prove that 
$(-\tilde\Delta+\frac{\partial}{\partial t})w(s,t)\ge 0$ in a weak sense.
We recall the definition of the conformal factors $f$ and $a$ resp.
in Example~\ref{zweienden}.
We compute for $\phi(s,t)\in C_0^\infty\left(M\times(0,T)\right), \phi\ge 0$
\begin{eqnarray*}
\lefteqn{
\int_M\int_{(0,T)}w(s,t)
\left(\left(-\tilde\Delta^*-\frac{\partial}{\partial
t}\right)\phi(s,t)\right)\, f^m\, dtdx}\\
&\ge& a(0)^{m-1}\int_0^T\left(\frac{\partial w}{\partial
s}(0-0,t)-\frac{\partial w}{\partial s}(0+0,t)\right)\phi(0,t)dt\ge
0.
\end{eqnarray*}

Moreover, since $\frac{\partial}{\partial s}w(s,t)<0$ for $s>0$, we see that
$$
w(s,t)\le w(0,t)=t^c
F\left( -c,1+\frac{1}{2\delta},-\frac{1}{4}t^{-1}\right)<2t^c\seq 0
$$ 
for $t\gg 1$.

Using the maximum principle stated in Lemma \ref{max4} yields the statement.
\end{proof}

\appendix\section{A general local parabolic maximum principle}

For the reader's convenience we shall outline the derivation
of a  local parabolic maximum principle, which is even more
general than we need it in the proof of Theorem~\ref{globalexistence}. 
Of particular interest is the dependence
of the estimation constants among others on
the elliptic operator and  the size and shape of the
domains. To a large extent, we follow \cite[Ch. VII]{Lieberman}.

Let $\Omega\subset{\bb R}^{n}\times {\bb R}$ be a domain. We denote the coordinates 
$X=(x,t)\in{\bb R}^n\times{\bb R} \left(\mbox{ resp. } Y=(y,s)\right)$.
In this section
we consider the operator 
$$
Lu:=-u_t+a^{ij}D_{ij}u+b^iD_iu+cu
$$
with real valued bounded measurable coefficients. Moreover we assume
the symmetric matrix  $(a^{ij})$ to be positive semidefinite.
We abbreviate ${\cal D}:=\det \left( a^{ij} \right) $ 
and ${\cal D}^*:={\cal D}^{\frac 1{n+1}}$. Furthermore, $\Lambda(X)$ denotes the
maximal and $\lambda(X)$ the minimal eigenvalue of $a^{ij}(X)$.
The function $u$ is considered in $u\in W^{2,1}_{n+1,loc}(\Omega)\cap C^0(\overline\Omega)$. 
As usual we denote by ${\cal P}\Omega$ the
parabolic boundary and by ${\cal B}\Omega$ the bottom of the domain $\Omega$.

We define the upper contact set $E(u)$ to be the set of all 
$X\in\Omega\setminus{\cal P}\Omega$ such that there exists 
$\xi\in{\bb R}^{n}$ such that
\begin{equation}u(X)+\xi(y-x)\ge u(Y)\end{equation}
for all $Y$ with $s\le t$. This implies $u_t\ge 0, -D^2u\ge 0$ on $E(u)$.

If $\Omega=B_R\times(0,T)$, we write $E^+(u)$ for the subset of $E(u)$ in which $u> 0$ and
\begin{equation}R|\xi|<u(X)-\xi\cdot x<\sup_\Omega\frac{u^+}2.\end{equation}

Similarly, we denote by $\Sigma(u)$ the set of all $\Xi=(\xi,h)\in{\bb R}^{n+1}$ such that 
\begin{equation}
R|\xi|<h<\sup_\Omega \frac{u^+}{2}.
\end{equation}

\vspace{4mm}
\noindent
First we quote  the global version of a maximum principle involving $L^p$-norms. 

\begin{proposition}{\bf\cite[Theorem 7.1]{Lieberman}}\,\ \ \Punkt \label{firstmax}Let $\Omega=B_R\times(0,T)$ and 
$u\in C^{2,1}(\Omega)\cap C^0(\overline\Omega)$ satisfying
$Lu\ge f$ with $c\le k$ in $\Omega$, where $k$ is a nonnegative constant. Then
$$\sup_\Omega u\le \exp(kT)\left(\sup_{{\cal P}\Omega}u^++c_1(n)B_0
R^{\frac n{n+1}}\left\|\frac f{{\cal D}^*}\right\|_{n+1,E^+(w)}\right),$$
with $B_0:=R^{-1}\|\frac b{{\cal D}^*}\|^{n+1}_{n+1,E^+(w)}+1$ and 
$w(x,t):=\exp(-kt)u-\sup_{{\cal P}\Omega}\left(\exp(-k\, . \,)u^+\right)$.
\end{proposition}

\vspace{4mm}
\noindent
Our goal is to prove the local counterpart of the preceding result. The crucial point
will be to estimate $\sup u$ by $\|f\|_{n+1}$ and 
 the weakest possible ``norm'' of $u$.
Since we will argue by means of  a scaling argument in the next proof, let us consider the degrees of 
the coefficients with respect to the two-parameter group ${\bb R}^2\cong\left(x\mapsto kx,t\mapsto lt\right)$. 
A simple calculation shows:
\begin{eqnarray*}\deg R&=&(-1,0)\\
\deg T&=&(0,-1)\\
\deg a^{ij}&=&(-2,1)\\
\deg b&=&(-1,1)\\
\deg c&=&(0,1)\\
\deg f&=&(0,1)
\end{eqnarray*}

For the following result, 
cf. \cite[Theorem 7.21]{Lieberman}.

\begin{proposition}[Local parabolic maximum principle]\Punkt\label{parmax}Let 
$\Omega=B_R\times(-T,0)$ and $u\in W^{2,1}_{n+1,loc}(\Omega)\cap C^0(\overline\Omega)$ satisfying
$Lu\ge f$. Assume further that 
$$\lambda\ge\lambda_0 >0,\quad \Lambda\le\Lambda_0,\quad |b|\le B,\quad c\le c_0.$$
Then for any $p>0$ and $0<\rho<1$ there exists $C$ depending only on $p,\rho$ and
$$(\lambda_0^{-n}TR^{-2})^{\frac 1{n+1}}(c_0R^2+BR+\Lambda_0+R^2T^{-1})$$ such that
$$\sup_{\rho\Omega} u\le C\left(|\Omega|^{-\frac 1p}\|u^+\|_p+(TR^{-1})^{\frac n{n+1}}\|f\|_{n+1}\right).$$
\end{proposition}

\begin{proof} By approximation, we may assume that 
$u\in C^{2,1}(\Omega)\cap C^0(\overline\Omega)$.
We note that both sides of the claimed inequality are invariant under the scaling 
$$x\mapsto lx, t\mapsto kt.$$ Hence it
suffices to prove the theorem for $\Omega=Q(1)=B_1\times[-1,0]$. 

For this purpose we define $\zeta:=(1-|x|^2)^+(1+t)^+$ and $\eta:=\zeta^q$ for $q>2$. 
We define the operator $P$ as principal part of $L$ by
$$Pv:=-v_t+a^{ij}D_{ij}v.$$
We will apply it to $v=\eta u$. This yields
$$Pv\ge \eta f-\eta(b^iD_iu+cu)+uP\eta+2a^{ij}D_iuD_j\eta.$$
We will calculate the terms separately. First we note that by Cauchy's inequality 
\begin{equation}
|Dv|\le\frac{v}{1-|x|}  \mbox{ on }E^+(v).\label{par0}
\end{equation}
{From} this it is easy to see that on $E^+(v)$
\begin{equation}
|Du|\le2(1+q)\frac v{\zeta\eta}.\label{par1}
\end{equation}

In order to compute $uP\eta$ -- again on $E^+ (v)$ -- 
we first note that $\eta_t\le q\frac \eta\zeta$. Next we use $(a^{ij})\ge 0$ to conclude
$$a^{ij}D_{ij}\eta\ge-2q\left(\mbox{tr\,}(a^{ij})\right)\frac \eta\zeta,$$
hence
\begin{equation}
uP\eta\ge -q(1+2\mbox{tr\,}a^{ij})\frac v\zeta\ge -q(1+2n\Lambda)\frac v\zeta.\label{par2}
\end{equation}

Finally, we have to compute $a^{ij}D_iuD_j\eta$. This splits up into the sum 
$\frac{a^{ij}}{\eta}D_ivD_j\eta-\frac{a^{ij}}{\eta^2}vD_i\eta D_j\eta$. 
For the first summand we obtain on $E^+(v)$ using (\ref{par0})
$$
\left| \frac{a^{ij}}{\eta}D_ivD_j\eta \right|\le 4q\Lambda\frac{v}{\zeta^2}.
$$
The second summand can be estimated
$$\left|\frac{a^{ij}}{\eta^2}vD_i\eta D_j\eta\right|\le\frac{\Lambda v|D\eta|^2}{\eta^2}\le 4q^2\Lambda\frac v{\zeta^2},$$
hence
\begin{equation}a^{ij}D_iuD_j\eta\ge -4q\Lambda(1+q)\frac v{\zeta^2}.\label{par3}\end{equation}
Adding up (\ref{par1}),(\ref{par2}) and (\ref{par3}) yields
\begin{eqnarray}Pv&\ge&\eta f-v\zeta^{-2}\left(c\zeta^2+\left(2(1+q)|b|+q(1+2n\Lambda)\right)\zeta
               +8q(1+2q)\Lambda\right)\nonumber\\
&\ge&\eta f-\tilde Cv\zeta^{-2}{\cal D}^*\label{par}
\end{eqnarray}
on $E^+(v)$, where $\tilde C$ can be chosen as
$$\tilde C:=2q(4(1+q)+n)\, \lambda_0^{-n/(n+1)}\, (c_0+B+\Lambda_0+1).$$
Note that the unique homogenization of $C$ to an element of degree $(0,0)$ (in $(k,l)$) gives the form mentioned in the
theorem.
Since $v=0$ on ${\cal P}\Omega$ and $P$ is an operator with $b=c=0$ we can apply Proposition \ref{firstmax} with $w=v$ and $k=0$.
This yields with $c_1=c_1 (n)$:
$$\sup_\Omega v\le c_1\left(\left\|\frac{f}{{\cal D}^*}\right\|_{n+1}+\tilde C\|v\zeta^{-2}\|_{n+1}\right)\le
c_1\tilde C(\|f\|_{n+1}+\|v\zeta^{-2}\|_{n+1}).$$

If $p>n+1$, we use H\"older's inequality, $v\zeta^{-2}\le u^+$ and $\eta\ge (1-\rho)^{2q}$ on $\rho\Omega$
to conclude the claim of the theorem. Here we may choose $q=2$.

If $p\le n+1$, we note that $v\zeta^{-2}=u^{\frac 2q}v^{1-\frac 2q}$. We choose $q=\frac{2(n+1)}p$ and compute 
$$\|v\zeta^{-2}\|\le (\sup_\Omega v)^{1-\frac p{n+1}}\|u^+\|_p^{\frac p{n+1}}\le \varepsilon\sup_\Omega v+c_2(n,p)\varepsilon^{1-\frac{n+1}p}
\|u^+\|_p$$
by Young's inequality.
Now we choose $\varepsilon:=(2c_1\tilde C)^{-1}$ and proceed like before. 
\end{proof}

\section{The energy differential inequality}\label{appendixb}

If $u:M\seq N$ satisfies the Hermitian harmonic system and $N$ has nonpositive
sectional curvature, in \cite{JostYau} an energy inequality
is mentioned which we are using several times.
This is why we sketch the proof here.

\begin{proposition}\Punkt If $N$ has nonpositive sectional curvature, then for every
relatively compact open set $\Omega\subset M$ there exists a constant
$C(\Omega)$ such that on $\Omega$
$$-\tilde\Delta e(u)\le C(\Omega)e(u).$$
\end{proposition}

\begin{proof}First we fix $x\in M$ and choose coordinates such that in $x$
resp. $f(x)$
$$\gamma_{\alpha\bar\beta}=\delta_{\alpha\bar\beta},\quad g_{ij}=\delta_{ij},\quad g_{ij,k}=0.$$
With these choices the left hand side becomes
\begin{eqnarray}\tilde\Delta e(u)&=&\sum_j \sum_\delta \Big(
\gamma^{\ia\qb}_{,\id\qd}u^j_{,\ia}u^j_{,\qb}\label{e1}\\
& & +\left[\gamma^{\ia\qb}_{,\id}u^j_{,\ia\qd}u^j_{,\qb}+\gamma^{\ia\qb}_{,\id}u^j_{,\ia}u^j_{,\qb\qd}\right]\label{e2}\\
& & +\left[\gamma^{\ia\qb}_{,\qd}u^j_{,\ia\id}u^j_{,\qb}+\gamma^{\ia\qb}_{,\qd}u^j_{,\ia}u^j_{,\qb\id}\right]\label{e3}\\
& & +\left[u^j_{,\ia\id}u^j_{,\qa\qd}+u^j_{,\ia\qd}u^j_{,\qa\id}\right]\label{e4}\\
& & +\left[u^j_{,\ia\id\qd}u^j_{,\qa}+u^j_{,\ia}u^j_{,\qa\id\qd}\right]\label{e5}\\
& & +(g_{ij,kl}\circ u)u^k_{,\id}u^l_{,\qd}u^i_{,\ia}u^j_{,\qa}\Big).\label{e6}
\end{eqnarray}

For (\ref{e1}) we obtain the inequality
$$\left|\sum_j \sum_\delta
\gamma^{\ia\qb}_{,\id\qd}u^j_{,\ia}u^j_{,\qb}\right|\le \Lambda e(u),$$
where $\Lambda$ can be estimated by bounds of terms of $\gamma_{\ia\qb}$ and their second derivatives.

In a similar way we can estimate (\ref{e2}) by
$$\left| \sum_j \sum_\delta
\gamma^{\ia\qb}_{,\id}u^j_{,\ia\qd}u^j_{,\qb}+\gamma^{\ia\qb}_{,\id}u^j_{,\ia}u^j_{,\qb\qd}
\right|\le
mC(\varepsilon|D^2u|^2+\frac{4m}{\varepsilon} e(u))$$
for any $\varepsilon>0$. The constant depends on $\gamma^{\ia\qb}$ and their
first derivatives.
The same holds true for (\ref{e3}).

It is not hard to see that (\ref{e4}) equals
$$
\sum_j \sum_\delta \left(
u^j_{,\ia\id}u^j_{,\qa\qd}+u^j_{,\ia\qd}u^j_{,\qa\id}\right)=\frac 18|D^2u|^2.
$$

For $(\ref{e5})+(\ref{e6})$ we have to use the Hermitian harmonic map system and the nonpositivity of the sectional curvature
to conclude that
\begin{eqnarray*}
\lefteqn{\sum_j \sum_\delta \left(
\left[u^j_{,\ia\id\qd}u^j_{,\qa}+u^j_{,\ia}u^j_{,\qa\id\qd}\right]+
(g_{ij,kl}\circ u)u^k_{,\id}u^l_{,\qd}u^j_{,\ia}u^j_{,\qa}\right)}\\
&=& -2\sum_j \sum_\delta
R_{ijkl}(u^i_{,\ia}u^j_{,\id}u^k_{,\qa}u^l_{,\qd}+u^i_{,\qa}u^j_{,\id}u^k_{,\ia}u^l_{,\qd})\ge 0.
\end{eqnarray*}

Putting all together yields
$$-\tilde\Delta e(u)\le \Lambda e(u) + 
2mC\left(\varepsilon|D^2u|^2+\frac{4m}{\varepsilon} e(u)\right)-\frac 18|D^2u|^2.$$
Choosing $\varepsilon\le (16mC)^{-1}$ yields the claimed inequality.
\end{proof}

\vspace{6mm}\noindent
{\bf Acknowledgment.} We are grateful to Wolf von Wahl (University of Bayreuth)
for his suggestion to investigate Hermitian-harmonic maps on
noncompact manifolds.

\bigskip\noindent
Fakult\"at f\"ur Mathematik,
Universit\"at Magdeburg,
Postfach 4120,
D-39016 Magdeburg, Germany\\
Current address of M.K.: University of California San Diego, Department of Mathematics, 9500 Gilman Drive, La Jolla, California,
92093-0112, USA
\end{document}